\documentclass[11pt,longtable,letterpaper,oneside,reqno]{amsart}
\usepackage{amsmath,amsfonts,amsthm,amssymb,amscd,mathtools,stmaryrd,url}
\usepackage{bbm}
\usepackage{epsfig}
\usepackage{graphicx}
\usepackage{latexsym}
\usepackage{longtable}
\usepackage{float}
\usepackage{hyperref}
\usepackage{color}

\DeclareFontFamily{OT1}{rsfs}{}
\DeclareFontShape{OT1}{rsfs}{n}{it}{<-> rsfs10}{}
\DeclareMathAlphabet{\mathscr}{OT1}{rsfs}{n}{it}

\newtheorem {defn }{Definition}

\numberwithin{equation}{section}
\usepackage[margin=1in]{geometry}
\usepackage{lipsum}

\begin{document}
\title{Remark on a paper by Wu}
\author[E. Hasanalizade]{Elchin Hasanalizade}
\address{Department of Mathematics and Computer Science, University of Lethbridge, 4401 University Drive, Lethbridge, Alberta, T1K 3M4 Canada}
\email{e.hasanalizade@uleth.ca}
\begin{abstract}
In this short note, we give an affirmative answer to Wu's conjecture on practical numbers, which was posed in [X.-H. Wu, {\it Special forms and the distribution of practical numbers}, Acta Math. Hungar., {\bf 160}(2020), 405-411].
\end{abstract}

\date{\today}
\maketitle

A positive integer $n$ is said to be {\it practical} if all positive integers $m<n$ can be written as a sum of distinct divisors of $n$. Stewart \cite{S} showed that an integer $n\ge2$, $n=p_1^{\alpha_1}\cdots p_k^{\alpha_k}$, with primes $p_1<\ldots<p_k$ and integers $\alpha_j\ge1$ is practical if and only if $p_1=2$ and 
\begin{align*}
p_j\le1+\sigma(p_1^{\alpha_1}\cdots p_{j-1}^{\alpha_{j-1}}) \ (2\le j\le k)
\end{align*}
where $\sigma(n)$ denotes the sum of the positive divisors of $n$.

Practical numbers have many properties similar to those of primes. Let $P(x)$ denote the number of practical numbers up to $x$. Saias \cite{Sa} using sieve methods proved a Chebyshev-type theorem: there are constants $c_1$ and $c_2$ such that 
\begin{align*}
c_1\frac{x}{\log{x}}<P(x)<c_2\frac{x}{\log{x}}.
\end{align*}
Margenstern \cite{Ma} conjectured that $P(x)\sim c\frac{x}{\log{x}}$ for a suitable constant $c$. This conjecture was confirmed by Weingartner \cite{W1} with the estimate 
\begin{align*}
P(x)=\frac{cx}{\log{x}}\bigg\{1+O\bigg(\frac{\log\log{x}}{\log{x}}\bigg)\bigg\}
\end{align*}
where $1.336073<c<1.336077$.

Some conjectures which are widely open for primes, have been proven for practical numbers. Melfi \cite{M} proved two Goldbach-type conjectures for practical numbers: i) every even positive integer is a sum of two practical numbers; ii) there exist infinitely many practical numbers $n$ such that $n-2$ and $n+2$ are also practical. Recently Wu \cite{Wu} has proved that for all integers $n\ge4$ there are at least two pactical numbers in the interval $(n^2,(n+1)^2)$. In the same paper, he stated the following conjecture.

{\bf Conjecture}. For all integers $k\ge1$, there exists an integer $N$ such that the interval $(n^2,(n+1)^2)$ contains at least $k$ practical numbers for all $n\ge N$.

In this note we answer affirmatively this conjecture.

\begin{proof}
Let $\mathcal{A}$ be the set of positive integers containing $n=1$ and all those $n\ge2$ with prime factorization $n=p_1^{\alpha_1}\cdots p_k^{\alpha_k}$, $p_1<p_2<\ldots<p_k$ which satisfy $p_1=2$ and 
\begin{align*}
p_i\le p_1^{\alpha_1}\cdots p_{i-1}^{\alpha_{i-1}} \ (2\le i\le k).
\end{align*}
Clearly the set $\mathcal{A}$ is a subset of the set of practical numbers. Weingartner \cite[Corollary 5]{W2} proved that for $x>x_0$ the interval $[x-x^{0.4872}, x]$ contains at least $x^{0.4872}(\log{x})^{-9.557}$ members of $\mathcal{A}$. Now let $k\ge1$ be a fixed integer. From Weingartner's theorem it follows that there exists $x_1=x_1(k)>x_0$ such that for $x>x_1$ the interval $[x-\sqrt{x}, x]$ contains at least $\frac{k}{2}+1$ members of $\mathcal{A}$. Let $N=\lfloor x_1\rfloor+1$. Then for $n\ge N$ there are at least $\frac{k}{2}+1$ members of $\mathcal{A}$ in each of the intervals 
\begin{align*}
\bigg[\bigg(n+\frac{1}{2}\bigg)^2-\bigg(n+\frac{1}{2}\bigg), \bigg(n+\frac{1}{2}\bigg)^2\bigg] \ \text{and} \ [(n+1)^2-(n+1), (n+1)^2].
\end{align*}
Since $\lceil n^2-\frac{1}{4}\rceil=n^2$ and $\lfloor n^2+n+\frac{1}{4}\rfloor=n^2+n$, the same conclusion is true for the intervals $[n^2, n^2+n]$ and $[n^2+n, (n+1)^2]$. Thus $(n^2, (n+1)^2)$ contains at least $k$ practical numbers since if $n\in\mathcal{A}$, then $n^2\in\mathcal{A}$ and out of three numbers $n^2, n^2+n, (n+1)^2$ at most two can be practical.
\end{proof}


\end{document}